\documentclass[10pt]{article}
\usepackage{amstext, amsmath, amsthm, amssymb}
\usepackage{amsmath}
\usepackage{geometry} 
\usepackage{color}

\numberwithin{equation}{section}
\usepackage[nottoc]{tocbibind}
\usepackage{hyperref}
\usepackage{makeidx}
\makeindex

\newtheorem{thm}{Theorem}[section]
\newtheorem{cor}[thm]{Corollary}
\newtheorem{lem}[thm]{Lemma}
\newtheorem{prop}[thm]{Proposition}
\newtheorem{defin}[thm]{Definition}
\newtheorem{rem}[thm]{Remark}
\newtheorem{example}[thm]{Example}

\usepackage{colortbl}

\def\liml{\lim\limits}

\def\suml{\sum\limits}
\def\intl{\int\limits}

\def\ker{\operatorname{Ker}}
\def\span{\operatorname{span}}

\title{Functional analysis approach to the Collatz conjecture
}

\author{Mikhail Neklyudov}

\AtEndDocument{\bigskip{\footnotesize
 (M.\,Neklyudov) \textsc{Instituto de Ci\^{e}ncias Exatas, Departamento de Matematica, UFAM, Manaus, CEP 69077-000, Brasil} \par
  \textit{E-mail address}: \texttt{misha.neklyudov@gmail.com}
}}

\date{}

\begin{document}
\maketitle

\begin{abstract}
We investigate the problems related to the Collatz map $T$ from the point of view of functional analysis. We associate with $T$ certain linear operator $\mathcal{T}$ and show that cycles and (hypothetical) diverging trajectory (generated by $T$) correspond to certain classes of fixed points of operator $\mathcal{T}$. 
Furthermore, we demonstrate connection between dynamical properties of operator $\mathcal{T}$ and map $T$. We prove that absence of nontrivial cycles of $T$ leads to hypercyclicity of operator $\mathcal{T}$. In the second part we show that the index of operator $Id-\mathcal{T}\in\mathcal{L}(H^2(D))$ gives upper estimate on the number of cycles of $T$. For the proof we consider the adjoint operator $\mathcal{F}=\mathcal{T}^*$
\[
\mathcal{F}: g\to g(z^2)+\frac{z^{-\frac{1}{3}}}{3}\left(g(z^{\frac{2}{3}})+e^{\frac{2\pi i}{3}}g(z^{\frac{2}{3}}e^{\frac{2\pi i}{3}})+e^{\frac{4\pi i}{3}}g(z^{\frac{2}{3}}e^{\frac{4\pi i}{3}})\right),
\]
first introduced by Berg, Meinardus in \cite{BM1994}, and show it does not have non-trivial fixed points in $H^2(D)$.

Moreover, we calculate resolvent of operator $\mathcal{F}$ and as an application deduce equation for the characteristic function of total stopping time $\sigma_{\infty}$.

%

 
Furthermore, we construct an invariant measure for  $\mathcal{T}$ in a slightly different setup, 
and investigate how the operator $\mathcal{T}$ acts on generalized arithmetic progressions.

\medskip

\medskip

\medskip

\end{abstract}

\section{Introduction}

Collatz conjecture has been approached from the multitude of points of view such as dynamical systems, Markov chains theory and many others (see the book \cite{Lagarias2010} for review and references). The idea to apply function theory and complex analysis to the Collatz conjecture first appeared in the papers of Berg, Meinardus \cite{BM1994} where they deduce functional equation for characteristic function of iterations of Collatz map and equivalent formulation of Collatz conjecture in terms of multiplicity of eigenvalue $1$ of certain linear operator $\mathcal{F}$ given by formula
\[
\mathcal{F}: g\to g(z^2)+\frac{z^{-\frac{1}{3}}}{3}\left(g(z^{\frac{2}{3}})+e^{\frac{2\pi i}{3}}g(z^{\frac{2}{3}}e^{\frac{2\pi i}{3}})+e^{\frac{4\pi i}{3}}g(z^{\frac{2}{3}}e^{\frac{4\pi i}{3}})\right).
\]
which can be characterized by 
\[
\mathcal{F}:\sum\limits_{n=0}^{\infty} a_n z^n\to \sum\limits_{n=0}^{\infty} a_{T(n)} z^n
\]
In this paper we will work mainly with operator $\mathcal{T}$ which can be characterised by
\[
\mathcal{T}:\sum\limits_{n=0}^{\infty} a_n z^n\to \sum\limits_{n=0}^{\infty} a_{n} z^{T(n)}
\]
The operators $\mathcal{F}$ and $\mathcal{T}$ are adjoint to each other in Hardy space $H^2(D)$, where $D$ is an open unit disk (Lemma \ref{lem:TFadjoint}). We will mainly concentrate on the connections between between dynamical properties of operator $\mathcal{T}$ and Collatz conjecture. The main result is an upper bound on number   of cycles of Collatz map $T$ by index of operator $Id-\mathcal{T}$ (Theorem \ref{thm:index}). In spite of duality between operators $\mathcal{F}$ and $\mathcal{T}$ they seems to relate to Collatz conjecture in a non-equivalent way. The reason is that in the natural set-up for the criterion proved by Berg, Meinardus \cite{BM1994}(Theorem $5$) operator $\mathcal{F}$ must be considered on the larger space of holomorphic functions $A(D)$ (with topology of uniform convergence on compact subsets) since non-trivial fixed points of $\mathcal{F}$ do not belong to Hardy space $H^2(D)$.


Let $T:\mathbb{Z}\to\mathbb{Z}$  be the Collatz map (reduced) defined by 
\[
T(n):=
\left\{
\begin{array}{rcl}
\frac{3n+1}{2} & , & n\in 2\mathbb{Z}+1\\
\frac{n}{2} &  , & n\in 2\mathbb{Z}
\end{array}
\right.
\]
and $H_{ber}^2(D)$ be a Bergman space on the open unit disk $D$. Define linear operator $\mathcal{T}: H_{ber}^2(D)\to H_{ber}^2(D)$ as follows
\[
\mathcal{T}f(z)= (Sf)(\sqrt{z})+\sqrt{z} (Af)(z^{\frac{3}{2}})
\]
where
\[
Sf(z):=\frac{f(z)+f(-z)}{2}, Af(z):=\frac{f(z)-f(-z)}{2} 
\]
are symmetric and antisymmetric parts of $f$.

Consequently, we could deduce that 
\[
\mathcal{T}(z^n):=z^{T(n)}, n\in \mathbb{Z}.
\]

It is easy to see that each cycle of the map $T$ corresponds to fixed point of the map $\mathcal{T}$ of polynomial form:
\begin{lem}\label{lem:Polynomialcycle}
If $(n_1,n_2,\ldots,n_k)$ is a cycle of map $T$  then 
\begin{equation}\label{eqn:Polynomialcycle}
\sum\limits_{i=1}^{k} z^{n_i}
\end{equation}
is a fixed point of operator $\mathcal{T}$. Furthermore, if $p$  is a fixed point $\mathcal{T}$ and $p$ is a polynomial then $p$ is of the form \eqref{eqn:Polynomialcycle} up to kernel of $\mathcal{T}$ (see \eqref{eqn:Kernel} below for the kernel of $\mathcal{T}$).
\end{lem}
\begin{proof}
It is trivial and omitted.
\end{proof}
There are several natural questions which appears in connection with Lemma \ref{lem:Polynomialcycle} and definition of operator $\mathcal{T}$. 
For instance, can we relate possible diverging trajectories corresponding to the map $T$ to fixed points of operator $\mathcal{T}$? Do there exists any other fixed points which does not correspond to cycles or diverging trajectories? Can we relate dynamic properties of operator $\mathcal{T}$ with dynamical properties of map $T$?

In this paper we give certain partial results to the questions above. We construct explicitly fixed point (actually countable number of them) of operator  $\mathcal{T}$ corresponding to (hypothetical) diverging trajectory of the map $T$ (Lemma \ref{lem:DivSequenceFPoint}). We construct other fixed points of operator $\mathcal{T}$ independent of cycles or diverging trajectories (Theorem \ref{thm:ExistenceFixedPoints}). The criteria to differ these classes of fixed points is an open problem. From the dynamical systems point of view we show that if map $T$ doesn't have nontrivial cycles then operator $\mathcal{T}$ is hypercyclic on certain natural Banach space (Theorem \ref{NocyclesHiper}). The question if $\mathcal{T}$ is chaotic (i.e. the set of periodic points is dense) is open. 

In the second part, we calculate resolvent of operator $\mathcal{F}$ (Theorem \ref{thm:FunctionalEquation}). The result of Berg, Meinardus \cite{BM1994} in this framework  is equivalent to calculation of  resolvent of operator $\mathcal{F}$ in one particular point. As an application of the exact form of resolvent we get an equation  
for the characteristic function of total stopping time. Furthermore, applying the fact that $\mathcal{F}$ is an expansive map we derive upper bound for the number of cycles of Collatz map through index of operator $Id-\mathcal{T}$ (as an operator in Hardy space $H^2(D)$).


Also we show that operator $\mathcal{T}$ has natural invariant measure in a slightly different setup (Theorem \ref{thm:InvariantMeasure}). In the end, we investigate action of operator $\mathcal{T}$ on generalized arithmetic progressions (section \ref{sec:ArithmeticProgressions}), which leads to another example of fixed point of $\mathcal{T}$ (Theorem \ref{thm:HigherFixedPoint}).

Our approach to the Collatz problem has certain natural limitations. Calculation of operator $\mathcal{T}$ on monomial $z^{\alpha},\alpha\in \mathbb{C}$ gives us
\[
\mathcal{T} (z^{\alpha})= z^{\frac{\alpha}{2}}\frac{1+e^{\alpha\pi i}}{2}+z^{\frac{3\alpha+1}{2}}\frac{1-e^{\alpha\pi i}}{2}
\]
and we can immediately see that if $\alpha$ is not an integer,i.e. if $\alpha$ is a rational with odd denominator, then $\mathcal{T} (z^{\alpha})$ will not be equal $z^{T(\alpha)}$ (where $T$ is an extension to the ring of rationals with odd denominators). This shows us that the operator $\mathcal{T}$ cannot be directly related with extension of  Collatz map $T$ to $2$-adic integers (\cite{MW84}) as well as  an extension of $T$ to the continuous map (\cite{Chamberland96}). In the same time it poses an interesting\footnote{possibly related to quantum mechanics. Notice that $\left|\frac{1+e^{\alpha\pi i}}{2}\right|^2+\left|\frac{1-e^{\alpha\pi i}}{2}\right|^2=1$ for $\alpha\in\mathbb{R}$ and, therefore, could be interpreted as probability amplitudes of some kind of deformed oscillator} problem of understanding of limiting behaviour of iterates $\{\mathcal{T}^k z^{\alpha}\}_{k\in\mathbb{N},\alpha\in\mathbb{C}}$ as $k\to\infty$.

{\bf Acknowledgements:} I would like to thank Slava Futorny, Wolfgang Bock and Sergei Nekrashevich  for constant support, enlightening discussions and encouragement!

\section{Connections between Collatz map $T$ and operator $\mathcal{T}$} 
In this section we will consider
\[
\mathcal{T}:H_{ber}^2(D)\slash X\mapsto H_{ber}^2(D)\slash  X
\]
where $X=\span\{1,z,z^2\}$--trivial invariant subspace of $\mathcal{T}$. It is easy to see that 
\[
||\mathcal{T}||_{\mathcal{L}(H_{ber}^2(D)\slash X)}\leq 2.
\]
Thus the spectrum of operator $\mathcal{T}$ belongs to the closed disk $2\overline{D}\subset\mathbb{C}$. We can notice that $0$ is eigenvalue of infinite multiplicity. Indeed, $\mathcal{T}(z^{2k+1}-z^{6k+4})=0, k\in \mathbb{N}\cup \{0\}$ and 
\begin{equation}\label{eqn:Kernel}
\ker \mathcal{T}=\span\{z^{2k+1}-z^{6k+4},k\in \mathbb{N}\cup \{0\}\}
\end{equation}
\begin{lem}
Operator $\mathcal{T}\in\mathcal{L}(H_{ber}^2(D)\slash X)$ is surjective and $3n+1$ conjecture implies the following statement
\begin{equation}\label{eqn:Nilpotency}
H_{ber}^2(D)\slash X=\overline{\bigcup\limits_{n=1}^{\infty}\ker(\mathcal{T}^n)} 
\end{equation}
\end{lem}
\begin{proof}
Surjectivity follows from
\[
\mathcal{T} z^{2k}=z^k, k\in\mathbb{Z}^+.
\]

Notice that $3n+1$ conjecture is equivalent to the fact that all monomials belong to $\bigcup\limits_{n=1}^{\infty}\ker(\mathcal{T}^n)$. Consequently,  the set of polynomials belongs to $\bigcup\limits_{n=1}^{\infty}\ker(\mathcal{T}^n)$ and the statement \eqref{eqn:Nilpotency} follows.
\end{proof}

\begin{thm}\label{NocyclesHiper}
If $T$ has no nontrivial cycles than $\mathcal{T}$ is hypercyclic.
\end{thm}
\begin{proof}
We will consider only the case when $3n+1$ conjecture is satisfied. If there exists diverging sequence $\{T^k(m)\}_{k=1}^{\infty}$  the proof is similar.

We will use criteria of hypercyclicity (Corollary 1.5 p. 235 in \cite{GodefroyShapiro1991}) with 
$T=\mathcal{T}$, $Sg(z):=g(z^2)$. 
\[
S(z^m):=z^{2m}, m\geq 3.
\]
Then assumption about $T$ is satisfied by \eqref{eqn:Nilpotency}, $T\circ S=Id$ and assumption about $S$ follows from elementary calculations:
\[
||S^n(z^k)||_{H_{ber}^2(D)\slash X}^2=||z^{2^n k}||_{H_{ber}^2(D)\slash X}^2=\frac{\pi}{2^n k+1}\to 0,n\to \infty
\]
\end{proof}

Let $A(D)$ be the space of analytic functions on the open disk $D$.
\begin{thm}\label{thm:ExistenceFixedPoints}
Any $\lambda\in \sqrt{2}D$ is eigenvalue of operator $\mathcal{T}:H_{ber}^2(D)\slash X\to H_{ber}^2(D)\slash X$ of infinite multiplicity. Furthermore, for any $\lambda\in \mathbb{C}\setminus \sqrt{2}D$  there exists infinite sequence $\{h_m(\lambda,\cdot)\}_{m=1}^{\infty}\subset A(D)$ such that
\[
\mathcal{T}h_m(\lambda,\cdot)=\lambda h_m(\lambda,\cdot).
\]
\end{thm}

\begin{proof}
Define
\begin{equation}\label{eqn:LacunaryF}
g(\lambda,z):=\sum\limits_{p=0}^{\infty}\lambda^p z^{2^p}\in A(D), f_m(\lambda, z):= g(\lambda,z^m), h_m(\lambda, \cdot):=f_{6m+4}(\lambda, \cdot)-f_{2m+1}(\lambda, \cdot), m\in\mathbb{N}.
\end{equation}
Then elementary calculation shows that
\begin{equation}\label{ActionTau_fm}
\mathcal{T} f_m=\lambda f_m+ z^{T(m)}, m\in\mathbb{N},
\end{equation}
and, consequently,
\[
\mathcal{T} h_m=\lambda h_m, m\in\mathbb{N}.
\]
Now it is enough to notice that $||h_m(\lambda, \cdot)||_{H_{ber}^2(D)\slash X}<\infty$ if $|\lambda|<\sqrt{2}$.
\end{proof}
\begin{lem}\label{lem:DivSequenceFPoint}
For every diverging sequence $\{T^k(m)\}_{k=1}^{\infty}$ of map $T$ we have fixed point of $\mathcal{T}$ of the following form:
\[
f_m(1,z)+\suml_{k=1}^{\infty}z^{T^k(m)}.
\]
\end{lem}
\begin{proof}
It immediately follows from formula \eqref{ActionTau_fm}.
\end{proof}
\begin{rem}\label{rem:DivSequenceFPoint}
We actually have infinite family of fixed points in Lemma \ref{lem:DivSequenceFPoint}. Indeed, instead of initial sequence $\{T^k(m)\}_{k=1}^{\infty}$ we could consider the truncated sequence  $\{T^k(m)\}_{k=n}^{\infty}$ for arbitrary $n\in\mathbb{N}$ and use the Lemma \ref{lem:DivSequenceFPoint} to construct fixed point for a new sequence.
\end{rem}
\begin{rem}
Thus we see that cycles and diverging sequences of $T$ correspond to fixed points of special form of operator $\mathcal{T}$. In the same time, operator $\mathcal{T}$ has other fixed points independently of existence of cycles and diverging sequences for $T$  as Theorem \ref{thm:ExistenceFixedPoints} shows. It would be interesting to understand how we can differ these classes of fixed points.
\end{rem}

\section{Invariant measure for the operator $\mathcal{T}$}\label{sec:Invmeasure}

The operator $\mathcal{T}$ can be written as follows
\[
\mathcal{T}=L(I+B)
\]
where
\begin{equation}\label{eqn:LBdef}
Lf(z):=\frac{1}{2}(f(\sqrt{z}))+f(-\sqrt{z})), Bf(z):=z f(z^3)
\end{equation}
and $I$ is an identity operator. Let us look at the behaviour of the operator $\mathcal{T}$
when we restrict on the functions on unit circle $|z|=1$. With parametrisation $z= e^{2\pi i\phi}$, $\phi\in [0,1)$ we have
\[
\mathcal{T} g(\phi)=\frac{1}{2}\left[g\left(\frac{\phi}{2}\right)+g\left(\frac{\phi+1}{2}\right)\right]+\frac{e^{\pi i\phi}}{2}\left[g\left(\frac{3\phi}{2}\right)-g\left(\frac{3\phi+1}{2}\right)\right],
\]
for $g\in L_{loc}^{2}(\mathbb{R},\mathbb{C})$, $g(\phi+1)=g(\phi),\phi\in\mathbb{R}$. As before, we have decomposition
\[
\mathcal{T}=L(I+B)
\]
where
\[
Lg(\phi):=\frac{1}{2}\left[g\left(\frac{\phi}{2}\right)+g\left(\frac{\phi+1}{2}\right)\right], Bg(\phi):=e^{2\pi i\phi} g(3\phi)
\]
Now we can see that $L$ is a standard transfer operator (\cite{LasMackey94},\cite{BezJorg2018}) with an invariant Lebesgue measure on the circle. Hence we can deduce that
\begin{thm}\label{thm:InvariantMeasure}
Lebesgue measure $dx$ on the circle (identified with interval $[0,1)$ as above) is invariant for the operator $\mathcal{T}:L^{2}(\mathbb{S}^1,\mathbb{C})\to L^{2}(\mathbb{S}^1,\mathbb{C})$ i.e.
\[
\int\limits_{\mathbb{S}^1}\mathcal{T}f\,dx=\int\limits_{\mathbb{S}^1}f\,dx, f\in L^{2}(\mathbb{S}^1,\mathbb{C}).
\]
\end{thm}
\begin{proof}
We have 
\[
\int\limits_{\mathbb{S}^1}\mathcal{T}f\,dx=
\int\limits_{\mathbb{S}^1}(f+Bf)\,dx
\]
by invariance of Lebesgue measure with respect to transfer operator $L$. Thus it remains to show that
\[
\int\limits_{\mathbb{S}^1}Bf\,dx=0.
\]
We have that
\[
\int\limits_{\mathbb{S}^1}Bf\,dx=\int\limits_0^1 e^{2\pi i\phi} f(3\phi)d\phi=\frac{1}{3}\int\limits_0^3 e^{\frac{2\pi i\psi}{3}} f(\psi)d\psi
\]
where $f$ is periodically extended beyond $[0,1)$.
Consequently, we have
\[
\int\limits_0^3 e^{\frac{2\pi i\psi}{3}} f(\psi)d\psi=
\left(\int\limits_0^1 +\int\limits_1^2+\int\limits_2^3 \right)e^{\frac{2\pi i\psi}{3}} f(\psi)d\psi=(1+\xi+\xi^2)\int\limits_0^1 e^{\frac{2\pi i\psi}{3}} f(\psi)d\psi,
\]
where $\xi=e^{\frac{2\pi i}{3}}$ and the last equality follows from periodicity of $f$. Now the result follows since $0=\xi^3-1=(\xi-1)(1+\xi+\xi^2)$ and $\xi\neq 1$.
\end{proof}

\section{Characterisation of total stopping time in terms of spectral properties of linear operator}

In this paragraph we will show that the index of operator $Id-\mathcal{T}\in\mathcal{L}(H^2(D))$ gives upper bound for the number of cycles of Collatz map $T$. For the proof we use expansiveness of operator $\mathcal{F}=\mathcal{T}^*\in\mathcal{L}(H^2(D))$  first introduced in Berg, Meinardis \cite{BM1994}.

They have calculated the value of resolvent of $\mathcal{F}$ at the point $\frac{z}{(1-z)^2}$ (p.7 formula 27 in \cite{BM1994}) and, as a result, they deduce an equation for characteristic function of Collatz iterations. Here the resolvent of $\mathcal{F}$ is calculated for certain class of analytic functions (Theorem \ref{thm:FResolvent}) and it allows us to find equation for characteristic function of total stopping time (Example \ref{ex:ex_1}). As a consequence, we show that  $\mathcal{F}$ has certain invariant subspace $Z$ explicitly written in \ref{def:Polynomials} and characterise Collatz conjecture in terms of long time dynamics of dynamical system generated by $\mathcal{F}$ on $Z$ (Corollary \ref{cor:Collatz}).


Let $H^2(D)$ be a Hardy space of holomorphic functions on open disk $D$ with norm
\[
||f||_{H^2(D)}^2=\sup\limits_{0\geq r<1}\left(\frac{1}{2\pi}\intl_0^{2\pi}|f(r e^{i\phi})|^2 d\phi\right)
\]
$\sigma_{\infty}:\mathbb{N}\cup\{0\}\to\mathbb{N}\cup\{\infty\}$ be a total stopping time (e.g. see \cite{Lagarias2010}, p.58) i.e minimal number (which is $\infty$ if it does not exist) such that $T^{\sigma_{\infty}(k)}(k)=1,k\in\mathbb{N}$. We will also put $\sigma_{\infty}(0):=0$.

We will need the following linear operator (first introduced in Berg, Meinardus \cite{BM1994}):
\begin{defin}
Define $\mathcal{F}: A(D)\to A(D-\{[-1,0]\})$ as follows
\begin{equation}\label{eqn:Operator_1}
\mathcal{F}(g)(z):=g(z^2)+\frac{z^{-\frac{1}{3}}}{3}\left(g(z^{\frac{2}{3}})+e^{\frac{2\pi i}{3}}g(z^{\frac{2}{3}}e^{\frac{2\pi i}{3}})+e^{\frac{4\pi i}{3}}g(z^{\frac{2}{3}}e^{\frac{4\pi i}{3}})\right)
\end{equation}
\end{defin}
\begin{prop}\label{prop:FBoundedness}
$\mathcal{F}\in \mathcal{L}(A(D))$ (where $A(D)$ is endowed with topology of uniform convergence on compact subsets), $\mathcal{F}$ is injective with closed range. Furthermore, $\mathcal{F}\in \mathcal{L}(H^2(D))$, $||\mathcal{F}||_{\mathcal{L}(H^2(D))}\leq \sqrt{2}$ and $\mathcal{F}$ is expansive in $H^2(D)$ i.e. $||f||_{H^2(D)}\leq ||\mathcal{F}f||_{H^2(D)},\, f\in H^2(D)$.
\end{prop}
\begin{proof}
\begin{itemize}
\item
Let us notice that
\begin{equation}\label{eqn:ActionFMon}
\mathcal{F}(z^n)=\left\{
\begin{array}{rcl}
z^{2n} & , & n=0,1 (\mbox{mod } 3)\\
z^{2n}+z^{\frac{2n-1}{3}} & , & n=2 (\mbox{mod } 3)
\end{array}
\right. , n\in\mathbb{Z}
\end{equation}
Thus $\mathcal{F}$ transfers set of polynomials $\mathcal{P}$ in itself. Then we have by definition of $\mathcal{F}$ that
\[
||\mathcal{F}f||_{C(\overline{B(0,r)})}\leq ||f||_{C(\overline{B(0,r^2)})}+ r^{-\frac{1}{3}}||f||_{C(\overline{B(0,r^{\frac{2}{3}})})},f\in \mathcal{P},0<r<1
\]
where we use the maximum principle for analytic functions to establish that the supremum of polynomial $\mathcal{F}f,f\in \mathcal{P}$ over $\overline{B(0,r)}$ is achieved on the  boundary of the disk $\overline{B(0,r)}$. Since $\mathcal{P}$ is dense in $A(D)$ and $A(D)$ is a Fr\'{e}chet space we have $\mathcal{F}\in \mathcal{L}(A(D))$.
\item
Injectivity of $\mathcal{F}$ follows from the following observation. By identity \ref{eqn:ActionFMon} we have
\begin{equation}\label{eqn:ActionFPol}
\mathcal{F}(\sum\limits_{n=0}^\infty a_n z^n)=
\sum\limits_{n=0}^\infty a_n z^{2n}+\sum\limits_{k=0}^{\infty} a_{3k+2} z^{2k+1}
\end{equation}
Thus symmetric part of $\mathcal{F}g$ is exactly $g(z^2)$. Consequently, if $\mathcal{F}g=0, g\in A(D)$ then $g(z^2)=0, z\in D$ and $g=0$.
\item
Let us show that $\mathcal{F}$ has closed range. Let $\{f_n\}_{n=1}^{\infty}\subset \mathcal{F}(A(D))$ i.e. $f_n=\mathcal{F}g_n, g_n\in A(D), n\in\mathbb{N}$ and $f_n\to f$ in $A(D)$. Then $\{f_n\}_{n=1}^{\infty}$ is a Cauchy sequence in $A(D)$ and we have
\[
||\mathcal{F}(g_n-g_m)||_{C(\overline{B(0,r)})}\to 0, n,m\to\infty, r<1.
\]
Consequently, applying that symmetric part of $\mathcal{F}g(z)$ is $g(z^2)$ we get
\[
||g_n-g_m||_{C(\overline{B(0,r^2)})}\to 0, n,m\to\infty, r<1.
\]
Since $r<1$ is arbitrary we deduce that the sequence $\{g_n\}_{n=1}^{\infty}$ is a Cauchy sequence in $A(D)$. Thus there exists $g\in A(D)$ $g_n\to g$ in $A(D)$ and 
$f_n=\mathcal{F}g_n\to \mathcal{F}g$ in $A(D)$. By uniqueness of the limit we conclude that $f=\mathcal{F}g,g\in A(D)$.

\item
By identity \eqref{eqn:ActionFPol}
\[
||\sum\limits_{n=0}^m a_n z^n||_{H^2(D)}^2\leq ||\mathcal{F}(\sum\limits_{n=0}^m a_n z^n)||_{H^2(D)}^2=
\sum\limits_{n=0}^m |a_n|^2+\sum\limits_{k=0}^{\left[\frac{m-2}{3}\right]} |a_{3k+2}|^2\leq 2\sum\limits_{n=0}^m |a_n|^2=2||\sum\limits_{n=0}^m a_n z^n||_{H^2(D)}^2
\]
and the result follows.
\end{itemize}
\end{proof}
\begin{lem}\label{lem:TFadjoint} 
For $\mathcal{T}\in \mathcal{L}(H^2(D))$ we have
\[
\mathcal{F}=\mathcal{T}^*
\]
\end{lem}
\begin{proof}
Let $g:=\sum\limits_{n=0}^{\infty} a_n z^n, f:=\sum\limits_{n=0}^{\infty} b_n z^n$. Then by identity \eqref{eqn:ActionFPol} we have
\[
(\mathcal{F}g,f)_{H^2(D)}=\sum\limits_{n=0}^{\infty} a_n \overline{b_{2n}}+a_{3n+2} \overline{b_{2n+1}}
\]
On the other side, we have 
\[
\mathcal{T}f=\sum\limits_{n=0}^{\infty} b_{2n} z^n+b_{2n+1} z^{3n+2},
\]
and, consequently,
\[
(g,\mathcal{T}f)_{H^2(D)}=\sum\limits_{n=0}^{\infty} a_n \overline{b_{2n}}+a_{3n+2} \overline{b_{2n+1}}=(\mathcal{F}g,f)_{H^2(D)}.
\]
\end{proof}
\begin{thm}\label{thm:index}
Number of cycles  of Collatz map $T:\mathbb{N}\to\mathbb{N}$ (including trivial one) is bounded above by index of operator $Id-\mathcal{T}\in \mathcal{L}(H^2(D))$.
\end{thm}
\begin{proof}
We have 
\[
ind (Id-\mathcal{T})=dim\ker(Id-\mathcal{T})-dim\ker(Id-\mathcal{T})^*=dim\ker(Id-\mathcal{T})-dim\ker(Id-\mathcal{F})
\]
Notice that $\mathcal{F}$ is an expansive operator by Proposition \ref{prop:FBoundedness}. We have that if $\mathcal{F}h=h, h\in H^2(D)$ then $||\mathcal{F}h||_{H^2(D)}^2=||h||_{H^2(D)}^2$ and we deduce from formula \eqref{eqn:ActionFPol} that $\mathcal{F}h(z)=h(z^2)$, which leads to $h$ being a constant function. Consequently, $dim\ker(Id-\mathcal{F})=1$. Now the result follows from Lemma \ref{lem:Polynomialcycle} and linear independence of polynomials corresponding to cycles (since cycles don't intersect).
\end{proof}
\begin{thm}\label{thm:FResolvent}
Let $l\in \mathbb{N}\cup\{0\}$, $\hat{\phi}=\sum\limits_{n=0}^{\infty}\phi(n)z^n\in A(D)$ such that 
\begin{equation}\label{eqn:phigrowth}
|\phi(n)|\leq C(n+1)^l,n\geq 0.
\end{equation}
Define 
\[
F^{\hat{\phi}}(z,w):=\sum\limits_{m,n=0}^{\infty} \phi(T^m(n))z^n w^m,z,w\in\mathbb{C}.
\]
Then $F^{\hat{\phi}}\in A(D\times \frac{2^l}{3^l}D)$ and
\begin{equation}
F^{\hat{\phi}}(z,w):=(I-w\mathcal{F})^{-1}\hat{\phi}.
\end{equation}
\end{thm} 
\begin{proof}
\begin{trivlist}
\item[(1)] Our assumption \eqref{eqn:phigrowth} on growth of $\phi$ together with an estimate 
\[
T^m(n)\leq (\frac{3}{2})^m (n+1)
\]
give us uniform convergence over compact subsets of $D\times \frac{2^l}{3^l}D$ of the series defining $F^{\hat{\phi}}$.
\item[(2)]
We have
\begin{eqnarray}
\mathcal{F}F^{\hat{\phi}} &= \sum\limits_{m,n=0}^{\infty} \phi(T^m(n))\mathcal{F}(z^n) w^m=\sum\limits_{m,n=0}^{\infty} \phi(T^m(n))z^{2n} w^m+\sum\limits_{m,k=0}^{\infty} \phi(T^m(3k+2))z^{2k+1} w^m\nonumber\\
&= \sum\limits_{m,n=0}^{\infty} \phi(T^{m+1}(2n))z^{2n} w^m+\sum\limits_{m,k=0}^{\infty} \phi(T^{m+1}(2k+1))z^{2k+1} w^m\nonumber\\
&= \sum\limits_{m,n=0}^{\infty} \phi(T^{m+1}(n))z^n w^m\nonumber\\
&= \frac{1}{w}(F^{\hat{\phi}}-\hat{\phi})\nonumber
\end{eqnarray}

\end{trivlist}
\end{proof}
We will need the following definition of characteristic function of total stopping time $\sigma_{\infty}$.
\begin{defin}
We call
\begin{equation}
g_{k,l}^{\lambda,\beta}(z):=\sum\limits_{m=0}^{\infty}\lambda^{\sigma_{\infty}(lm+k)}\beta^m z^{lm+k}, \lambda\in D\cup\{\ 1\},z\in D, k,l\in \mathbb{N},\beta\in (0,1]
\end{equation}
characteristic function of $\sigma_{\infty}$ with parameters $(k,l,\beta)$ and, corr.,
\[
\widetilde{g}_{k,l}^{\lambda}(z):=\sum\limits_{m=0}^{\infty}\lambda^{\sigma_{\infty}(lm+k)} z^m.
\]
reduced characteristic function of $\sigma_{\infty}$ with parameters $(k,l)$ (In this paragraph we use convention that $\lambda^{\infty}=0$ if $\lambda\in D$ and $1^{\infty}=1$). 
\end{defin}
\begin{rem}
Notice that $g_{k,l}^{1,\beta}=\frac{z^k}{1-\beta z^l}$. Thus, in certain sense, $g_{k,l}^{1,\beta}$ represents "smoothed out" (when $\beta\in (0,1)$) arithmetic progression $\{k+lp,p\in \mathbb{N}\}$.
\end{rem}
\begin{example}\label{ex:ex_1}
Let us consider
\[
\phi(n)=\delta_{n,1}
\]
in Theorem \ref{thm:FResolvent}.
Then $\hat{\phi}(z)=z$, $F^{\hat{\phi}}\in A(D\times D)$ and 
\[
F^{\hat{\phi}}=\sum\limits_{n=1}^{\infty}\left(\sum\limits_{m=0}^{\infty}\delta_{T^m(n),1}w^m\right)z^n=
\frac{1}{1-w^2}\sum\limits_{n=1}^{\infty} w^{\sigma_{\infty}(n)}z^n
\]
\end{example}
Consequently, from Theorem \ref{thm:FResolvent} we get
\begin{thm}\label{thm:FunctionalEquation}
Let $\lambda\in D-\{0\}$. Then 
 $\widetilde{g}_{0,1}^{\lambda}=g_{0,l}^{\lambda,1}$ satisfies the following identity
\begin{equation}\label{eqn:Operator_GenEig}
\mathcal{F}(\widetilde{g}_{0,1}^{\lambda}-1)=\frac{\widetilde{g}_{0,1}^{\lambda}-1}{\lambda}+ \left(\lambda-\frac{1}{\lambda}\right)z.
\end{equation} 
\end{thm}
The next example has been first proved in Berg, Meinardus (\cite{BM1994},p.7, formula 27). 
\begin{example}\label{ex:ex_2}
Let us consider
\[
\phi(n)=n
\]
in Theorem \ref{thm:FResolvent}. Then $\hat{\phi}(z)=\frac{z}{(1-z)^2}$, $F^{\hat{\phi}}\in A(D\times \frac{2}{3} D)$ and 
\[
F^{\hat{\phi}}=\sum\limits_{m,n=0}^{\infty} T^m(n) z^n w^m
\]
Theorem \ref{thm:FResolvent} implies
\[
F^{\hat{\phi}}-w\mathcal{F} F^{\hat{\phi}}=\frac{z}{(1-z)^2}
\]
\end{example}
\begin{rem}
At first glance it seems that the correction term $\left(\lambda-\frac{1}{\lambda}\right)z$ in formula \eqref{eqn:Operator_GenEig} can be deleted by considering summation over some subset of natural numbers in the definition of $\widetilde{g}_{0,1}^{\lambda}$. This is not the case, since existence of eigenvectors with eigenvalue $\frac{1}{\lambda},\lambda\in D-\{0\}$ would contradict boundedness of operator $\mathcal{F}$ in $A(D)$ (Proposition \ref{prop:FBoundedness}) when $|\lambda|$ is small enough.
\end{rem}
\begin{defin}\label{def:Polynomials}
Define family of polynomials
\begin{equation}
Pol_k(z)=\sum\limits_{n\in\mathbb{N}:\sigma_{\infty}(n)=k}z^n,k\in\mathbb{N}\cup\{0\}.
\end{equation}
Notice that the family of polynomials $\{Pol_k\}_{k=0}^{\infty}$ is linearly independent in $A(D)$. Denote a space
\[
Z:=\overline{\span\{Pol_k,k\in\mathbb{N}\cup\{0\}\}}\subset A(D),
\]
where closure is in the topology of $A(D)$. 
\end{defin}
\begin{cor}
$Z\subset A(D)$ is an invariant subspace for operator $\mathcal{F}$. Furthermore, $\mathcal{F}$ in the basis $\{Pol_k\}_{k=0}^{\infty}$ has the following representation:
\begin{equation}\label{eqn:FRepresentation}
\mathcal{F}=\left(
\begin{array}{cccccc}
0 & 1 & 0 & 0 & 0 & \ldots\\
1 & 0 & 1 & 0 & 0 & \ldots\\
0 & 0 & 0 & 1 & 0 & \ldots\\
0 & 0 & 0 & 0 & 1 & \ldots\\
0 & 0 & 0 & 0 & 0 & \ldots\\
\ldots & \ldots & \ldots & \ldots & \ldots & \ldots\\
\end{array}
\right) ,
\end{equation}
and $\mathcal{F}$ has simple eigenvalue $1$ (as an operator in $\mathcal{L}(Z)$) with eigenvector $\psi:=\sum\limits_{k=0}^{\infty} Pol_k=\sum\limits_{m\in\mathbb{N},\sigma_{\infty}(m)<\infty} z^m$.
\end{cor}
\begin{proof}
Combining identity
\[
\widetilde{g}_{0,1}^{\lambda}(z)=1+\sum\limits_{k=0}^{\infty}\lambda^k Pol_k(z), z\in D
\]
with equality \eqref{eqn:Operator_GenEig} we immediately get formula \eqref{eqn:FRepresentation}. Invariance of subspace and simplicity of eigenvalue $1$ with eigenvector $\psi$ follows. 
\end{proof}
The equivalence of Collatz conjecture and part (b) of the next corollary was first proved by L. Berg, G. Meinardus\cite{BM1994}. 

\begin{cor}\label{cor:Collatz}
The following are equivalent:
\begin{itemize}
\item[(a)] Collatz conjecture 

\item[(b)] The dimension of the vector space of eigenvectors of eigenvalue $1$ of operator $\mathcal{F}\in \mathcal{L}(A(D))$ is $2$.

\item[(c)] Eigenvector $\frac{z}{1-z}\in Z$.

\item[(d)] 

$\lim\limits_{n\to\infty}\mathcal{F}^n(z+z^2)=\frac{z}{1-z}$, where limit is in topology $A(D)$.

\end{itemize}
\end{cor}
\begin{proof}
\begin{itemize}
\item[(a)](a)$\Longrightarrow$(b). Suppose that Collatz conjecture is true. Let us consider  equation
\[
\mathcal{F}h=h,h\in A(D).
\]
We can put $h=\suml_{n=0}^{\infty} a_n z^n$ and rewrite it as follows
\[
\suml_{n=0}^{\infty} a_n z^n=\suml_{n=0}^{\infty} a_n z^{2n}+\suml_{k=0}^{\infty} a_{3k+2} z^{2k+1}.
\]
Then we have $a_{n}=a_{T(n)}, n\in \mathbb{N}\cup\{0\}$. Collatz conjecture implies that $a_n=a_1,n\in\mathbb{N}$.  Consequently, we have $h=a_0+a_1\frac{z}{1-z}$ and the result follows.
\item[(b)](b)$\Longrightarrow$(a). 
The operator $\mathcal{F}$ has at least two linearly independent fixed points: $\{1, \frac{z}{1-z}\}$.
Violation of Collatz conjecture implies existence of another fixed point $\psi=\sum\limits_{m\in\mathbb{N},\sigma_{\infty}(m)<\infty} z^m\in A(D)$ of operator $\mathcal{F}$ .
\item[(c)](a)$\Leftrightarrow$(c). 
If $\frac{z}{1-z}\in Z$ then by simplicity of eigenvalue $1$ of operator $\mathcal{F}\in\mathcal{L}(Z)$ we have $\psi=\frac{z}{1-z}$ and the conjecture follows. Conjecture implies $\frac{z}{1-z}=\psi\in Z$.
\item[(d)](c)$\Leftrightarrow$(d). 
It follows immediately from formula \eqref{eqn:FRepresentation} that 
\[
\mathcal{F}^n(Pol_0+Pol_1)=\sum\limits_{l=0}^{n-1} Pol_l
\]
Consequently, $\psi=\lim\limits_{n\to\infty}\mathcal{F}^n(z+z^2)$ and the result follows.

\end{itemize}
\end{proof}
\begin{rem}
Taking the limit $\lambda\to -1$ in identity \eqref{eqn:Operator_GenEig} we can get eigenvector
\[
\sum\limits_{m\in\mathbb{N},\sigma_{\infty}(m)<\infty} (-1)^{\sigma_{\infty}(m)}z^m\in A(D),
\]
of operator $\mathcal{F}$ with eigenvalue $-1$.
\end{rem}
\begin{rem}
It follows from \eqref{eqn:ActionFMon} that restriction of $\mathcal{F}$ on polynomials defines almost inverse (in a sense that both $\mathcal{F}\circ \mathcal{T}-I$, $\mathcal{T}\circ \mathcal{F}-I$ are degenerate operators) of operator $\mathcal{T}$.  This should be compared to the following well known reformulation of Collatz conjecture. 

Let $\Psi: 2^{\mathbb{N}}\to 2^{\mathbb{N}}$ be a $\sigma$-additive map defined by its values on one-element subsets of $\mathbb{N}$ as follows
\[
\Psi(\{n\})=\left\{
\begin{array}{rcl}
\{2n\} & , & n=0,1 (\mbox{mod } 3)\\
\{2n,\frac{2n-1}{3}\} & , & n=2 (\mbox{mod } 3)
\end{array}
\right. , n\in\mathbb{N}.
\]
Notice that 
\[
T\circ\Psi(\{n\})=\{n\}
\]
(as sets) i.e. $\Psi$ is a right inverse of $T$.
Then 
it follows that
\[
\Psi (\mathbb{N})=\mathbb{N}.
\]
At the same time, 
the set
\[
\mathcal{A}=\{n\in\mathbb{N}|\sigma_{\infty}(n)<\infty\}
\]
is also a fixed point of $\Psi$. Thus, Collatz conjecture is equivalent to uniqueness of fixed point of $\Psi$.
\end{rem}
\begin{cor}
For $q\in (0,\frac{1}{2})$ we have
\[
\sum\limits_{n=2}^{\infty} q^{\sigma_{\infty}(n)}\leq \frac{(2-q)q}{1-2 q}.
\]
\end{cor}
\begin{proof}
Let $h_{\lambda}=\widetilde{g}_{0,1}^{\lambda}-1$. By Theorem \ref{thm:FunctionalEquation} we have
\[
\mathcal{F}h_{\lambda}=\frac{1}{\lambda}h_{\lambda}+(\lambda-\frac{1}{\lambda})z
\]
If $\lambda$ satisfies $|\lambda|^2<\frac{1}{2}$ we apply   Lemma \ref{prop:FBoundedness} to deduce
\[
||\frac{1}{\lambda}h_{\lambda}+(\lambda-\frac{1}{\lambda})z||_{H^2(D)}^2=||\mathcal{F}h_{\lambda}||_{H^2(D)}^2\leq 2||h_{\lambda}||_{H^2(D)}^2.
\]
Consequently, we have
\[
\frac{1}{|\lambda|^2}\left(\sum\limits_{n=2}^{\infty}|\lambda|^{2\sigma_{\infty}(n)}+|\lambda|^4\right)\leq 2(1+\sum\limits_{n=2}^{\infty}|\lambda|^{2\sigma_{\infty}(n)}),
\]
and we can conclude that for $|\lambda|^2<\frac{1}{2}$
\[
\sum\limits_{n=2}^{\infty}|\lambda|^{2\sigma_{\infty}(n)}\leq \frac{|\lambda|^2(2-|\lambda|^2)}{1-2|\lambda|^2}.
\]
Denoting $q=|\lambda|^2$ we get the result.
\end{proof}


\section{Action of operator $\mathcal{T}$ on generalized arithmetic progressions}\label{sec:ArithmeticProgressions}

Let us look at how the operator $\mathcal{T}$ acts on functions $g_{k,l}^{\lambda,\beta}$ representing generalization of  arithmetic progression $\{k+ln\}_{n\in\mathbb{N}}$ (in some sense "weighted" with parameters $\lambda,\beta$). We have that
\begin{equation}\label{eqn:ActionGenProgrArithmetica}
\mathcal{T}g_{k,l}^{\lambda,\beta}=\lambda
\left\{
\begin{array}{rcl}
g_{\frac{k}{2},\frac{l}{2}}^{\lambda,\beta} & , & k,l\in 2\mathbb{Z}\\
g_{\frac{3k+1}{2},\frac{3l}{2}}^{\lambda,\beta} & , & k\in 2\mathbb{Z}+1,l\in 2\mathbb{Z}\\
\beta g_{\frac{3(k+l)+1}{2},3l}^{\lambda,\beta^2}+g_{\frac{k}{2},l}^{\lambda,\beta^2} & , & k\in 2\mathbb{Z},l\in 2\mathbb{Z}+1\\
\beta g_{\frac{k+l}{2},l}^{\lambda,\beta^2}+g_{\frac{3k+1}{2},3l}^{\lambda,\beta^2}& , & k,l\in 2\mathbb{Z}+1
\end{array}
\right.
\end{equation}
Natural candidate for fixed point of operator $\mathcal{T}$ will be the limit $\liml_{n\to\infty}\mathcal{T}^n g_{k,1}^{\lambda,\beta}$ since informally speaking it should be invariant w.r.t. $\mathcal{T}$. Now we will argue that the value of $\lambda$ should be $\frac{1}{2}$ to get nontrivial fixed point:
\begin{lem} We have for $0<\beta<1$ and $0<\lambda<\frac{1}{2}$ that  $\liml_{n\to\infty}\mathcal{T}^n g_{0,1}^{\lambda,\beta}=0$.
\end{lem}
\begin{proof}
It is easy to see that $\mathcal{T}^n g_{0,1}^{\lambda,\beta}$ has the form
\[
\mathcal{T}^n g_{0,1}^{\lambda,\beta}=\lambda^n\sum\limits_{i=1}^{2^n}g_{m_i,k_i}^{\lambda,\beta}
\]
where $\{m_i,k_i\}_{i=1}^{\infty}$ are some sequences of natural numbers. Now for $\lambda<\frac{1}{2}$ we have
\[
\liml_{n\to\infty}\mathcal{T}^n g_{0,1}^{\lambda,\beta}=0
\]
by uniform bound 
\[
|g_{m_i,k_i}^{\lambda,\beta}|\leq\frac{|\lambda|}{1-\beta z}.
\]
\end{proof}
\begin{rem}
Thus the natural candidate for the fixed point of $\mathcal{T}$ will be the limit $\liml_{n\to\infty}\mathcal{T}^n g_{0,1}^{\lambda_0,\beta}, \lambda_0=\frac{1}{2}$. It is an open problem to calculate this limit or, generally, calculate $\liml_{n\to\infty}\mathcal{T}^n g_{k,l}^{\lambda_0,\beta}, \lambda_0=\frac{1}{2}$, $k,l\in\mathbb{N}$. 
\end{rem}
From now on we will assume that $\lambda=\beta=1$. In this case we have $g_{k,l}^{1,1}=\frac{z^k}{1-z^l}$ i.e. $g_{k,l}^{1,1}$ corresponds precisely $\{k+ln\}_{n\in\mathbb{N}}$. Furthermore, it is easy to notice that the limit $\liml_{n\to\infty}\mathcal{T}^n g_{k,l}^{1,1}$ (natural candidate for fixed point of $\mathcal{T}$) does not exist for many progressions. For instance, if $l=1,k=0$ then we have repetition of the term $\frac{z^2}{1-z^3}$ in the sequence $\{\mathcal{T}^n g_{0,1}^{1,1}\}_{n=1}^{\infty}$ which leads to divergence. In this section we will suggest a way which will allow us to avoid the divergence and lead us to the construction of certain class of fixed points of operator $\mathcal{T}$ different from the construction of Theorem \ref{thm:ExistenceFixedPoints}.

The principal idea is to consider instead of arithmetic sequence $\{l+mn\}_{n\in\mathbb{N}}$ its subsequence $\{l+mn\}_{n\in\mathbb{K}_k}$, where $\mathbb{K}_k$ is a subset of natural numbers, which have in binary numeral system fixed number $k$ of ones  i.e. sums of $k$  different powers of $2$.

Fix $k\in \mathbb{N}$. To define the corresponding function let us look at the power $k$ of lacunary function $g(1,z)$ defined in \eqref{eqn:LacunaryF}. It is easy to see that
\[
\frac{g^k}{k!}=\suml_{0\leq p_1\leq p_2\leq\ldots\leq p_k }z^{2^{p_1}+2^{p_2}+\ldots+2^{p_k}}
\]
Consequently there exists polynomial $P_k(z)$ of degree $k$ with highest order term $\frac{z^k}{k!}$ such that
\[
P_k(g)=\suml_{0\leq p_1< p_2<\ldots< p_k }z^{2^{p_1}+2^{p_2}+\ldots+2^{p_k}},
\]
and we define the function corresponding to the  subsequence $\{l+mn\}_{n\in\mathbb{K}_k}$ as follows
\[
\Psi_{l,m}^k(z):=z^l P_k(g)(1,z^m).
\]
Also we put 
\[
P_0:=1, \Psi_{l,m}^0:=z^l.
\]
\begin{lem}\label{lem:ActionGenProgrArithmeticaSub}
\begin{equation}\label{eqn:ActionGenProgrArithmeticaSub}
\mathcal{T} \Psi_{l,m}^k=
\left\{
\begin{array}{rcl}
\Psi_{\frac{l}{2},\frac{m}{2}}^{k} & , & l,m\in 2\mathbb{Z}\\
\Psi_{\frac{3l+1}{2},\frac{3m}{2}}^{k} & , & l\in 2\mathbb{Z}+1,m\in 2\mathbb{Z}\\
\Psi_{\frac{l}{2},m}^{k}+\Psi_{\frac{3(l+m)+1}{2},3m}^{k-1} & , & l\in 2\mathbb{Z},m\in 2\mathbb{Z}+1\\
\Psi_{\frac{3l+1}{2},3m}^{k}+\Psi_{\frac{l+m}{2},m}^{k-1} & , & l,m\in 2\mathbb{Z}+1 
\end{array}
\right.
\end{equation}
\end{lem}

Comparing action of $\mathcal{T}$ on the arithmetic progression (formula \eqref{eqn:ActionGenProgrArithmetica} with $\lambda=\beta=1$) with action on its subsequence (formula \eqref{eqn:ActionGenProgrArithmeticaSub}) we notice that for the latter when we have a split the splitting function has parameter $k$ (which describes number of ones in the corr. subsequence) decreased to $k-1$. 
Based on this observation we can construct fixed point of $\mathcal{T}$ different from the ones constructed in Theorem \ref{thm:ExistenceFixedPoints}.

\begin{thm}\label{thm:HigherFixedPoint}
Function
\begin{equation}
FP_2:=\frac{g^2}{2}+z-g(z)+\suml_{k=1}^{\infty}\left[(z+z^2)g(z^{3^k})-g(z^{1+3^k})-g(z^{2+3^k})\right]\in A(D)
\end{equation}
is a fixed point of $\mathcal{T}$.
\end{thm}
\begin{proof}
Uniform convergence of the series inside of the disk follows from the following estimates
\begin{equation}
\left|\sum\limits_{k=1}^{\infty}g(z^{3^k})\right|\leq
\sum\limits_{k=1,l=0}^{\infty} |z|^{3^k2^l}\leq \sum\limits_{n=0}^{\infty}|z|^n\leq \frac{1}{1-|z|},
\end{equation}
\begin{eqnarray}
&\left|\sum\limits_{k=1}^{\infty}g(z^{1+3^k})\right|=
\left|\sum\limits_{k=1,l=0}^{\infty} z^{2^l+3^k2^l}\right|\nonumber\\
&\leq \sum\limits_{k=1}^{\infty}\suml_{l=0}^{\infty}|z|^{2^l+3^k2^l} \leq \sum\limits_{k=1}^{\infty}\suml_{l=0}^{\infty}|z|^{2^{l+1}}\suml_{l=0}^{\infty}|z|^{3^k2^{l+1}}\nonumber\\
&\leq 
\suml_{l=1}^{\infty}|z|^{2^l}\sum\limits_{k=1,l=1}^{\infty}|z|^{3^k2^l} \leq \frac{|z|^8}{(1-|z|)^2}
\end{eqnarray}
To show that $FP_2$ is a fixed point notice that by Lemma \ref{lem:ActionGenProgrArithmeticaSub} we have
\begin{equation}\label{eqn:Calc_1}
\mathcal{T}\Psi_{0,1}^2=\Psi_{0,1}^2+\Psi_{2,3}^1,
\end{equation}
\begin{equation}\label{eqn:Calc_2k}
\mathcal{T}\Psi_{1,3^k}^1=\Psi_{2,3^{k+1}}^1+ z^{T(1+3^k)},k\in\mathbb{N}, 
\end{equation}
\begin{equation}\label{eqn:Calc_3k}
\mathcal{T}\Psi_{2,3^k}^1=\Psi_{1,3^k}^1+ z^{T(2+3^k)}, ,k\in\mathbb{N}.
\end{equation}
Fix $M\in \mathbb{N}$. Taking sum over $k$ from $1$ to $M$ in \eqref{eqn:Calc_2k}, \eqref{eqn:Calc_3k} and adding \eqref{eqn:Calc_1} leads to
\begin{equation}\label{eqn:Calc_4}
\mathcal{T}\left[\Psi_{0,1}^2+\suml_{k=1}^{M}(\Psi_{1,3^k}^1+\Psi_{2,3^k}^1)\right]=\left[\Psi_{0,1}^2+\suml_{k=1}^{M}(\Psi_{1,3^k}^1+\Psi_{2,3^k}^1)
+\Psi_{2,3^{M+1}}^1\right]+\suml_{k=1}^{M}z^{T(1+3^k)}+z^{T(2+3^k)}
\end{equation}
On the other side by formula \eqref{ActionTau_fm} (with $\lambda=1$) we have
\begin{equation}\label{eqn:Calc_5}
\mathcal{T}\suml_{k=1}^{M}(f_{1+3^k}+f_{2+3^k})=\suml_{k=1}^{M}(f_{1+3^k}+f_{2+3^k})+\suml_{k=1}^{M}z^{T(1+3^k)}+z^{T(2+3^k)}.
\end{equation}
Subtracting from \eqref{eqn:Calc_4} formula \eqref{eqn:Calc_5} and taking the limit $M\to\infty$ we conclude the proof.
\end{proof}
\begin{rem}
We can observe that fixed point constructed in this paragraph as well as fixed points constructed in the Theorem \ref{thm:ExistenceFixedPoints} are analytic functions inside of the disk while fixed points corresponding to the cycle (polynomial) are analytic functions on the whole plane. Thus it would be natural to ask if there exist fixed point of operator $\mathcal{T}$ analytic on the whole plane different from $z+z^2$. Negative answer would imply no non-trivial cycles conjecture. Another natural question would be if we can find analytic function on $\mathbb{C}$ in the linear space generated by fixed points (constructed in the remark \ref{rem:DivSequenceFPoint}) corresponding to (hypothetical) diverging trajectory. 
\end{rem}

\end{document}